\documentstyle[aps,12pt]{revtex}
\begin{document}
\textwidth=14cm
\textheight=22cm
\title{Two-Parameter Quantum Groups and Quantum Planes }
\author{Sunggoo Cho and   Sang-jun Kang,}
\address{ Department of Physics,  Semyung University,
  Chechon,  Chungbuk 390-230,   Korea   }
\author{ Chung-hum Kim }
\address{ Department of Physics, Sun Moon University,
Chonan, Chungnam 330-150, Korea  }
\author{Kwang Sung Park}
\address{Department of Mathematics, Keimyung University,
	 Taegu 705-701, Korea}
\date{\today }
\maketitle
\draft
\hspace{2cm}
\begin{abstract}
Usually the generators of a quantum group are assumed to be commutative
with the noncommuting coordinates of a quantum plane.
We have relaxed the assumption and investigated its consequences.
Not only does a two-parameter quantum group arise naturally,
but also
the formulation leads us  to many probable
quantum planes associated with a quantum group.
Several examples are presented.

\end{abstract}
\newpage
 \section{Introduction}

In  recent years, the concept of a quantum group has extensively emerged
in the physical and mathematical literature \cite{ref1,ref2,ref3}.
  Quantum groups are nontrivial
generalizations
of ordinary Lie groups.  Such generalizations are made in the framework of
Hopf algebras \cite{ref4,ref5,ref6}.
 A Hopf algebra is an algebra together with
operations called the comultiplication, counit and antipode, which reflect
the group structure.
A quantum group is a non-commutative Hopf algebra consistent with these
costructures. Usually, quantum groups are introduced as deformations of
commutative Hopf  algebras in the sense that they become commutative Hopf
algebra as some parameters go to particular values\cite{ref7,ref8}.
Probably  the most
studied case of a quantum group is $ GL_q (2) $ whose element
$ T = \left( \begin{array}{cl} A & B \\ C & D \end{array} \right) $
satisfies the following nontrivial commutation relations:
\begin{eqnarray}
  AB &=& q BA, \hspace{2cm}  AC = q CA,    \nonumber \\
  BD &=& q DB, \hspace{2cm}  CD = q DC,     \label{eq:1}\\
  BC &=& CB, \hspace{2cm}  AD - DA = (q - q^{-1} ) BC \, . \nonumber
\end{eqnarray}
On the other hand, quantum spaces or quantum planes may be introduced
 as representation spaces
of quantum groups\cite{ref1,ref3,ref14}.

 Corresponding to the quantum group $GL_q (2)$,
 Manin \cite{ref1}
has defined a quantum space as  one
generated by two
noncommuting coordinates $x , \; y$ obeying
\begin{equation}
x \, y \,= \, q \, yx \; \; \; (q \, \neq \, 0 \, , \; 1)\, . \label{eq:2}
\end{equation}
Then the quantum group $GL_q (2)$ becomes a symmetry group of the quantum
plane. In fact, the points $(x^{\prime}  , \; y^{\prime} \, ) $ and
$(x^{\prime\prime}  , \; y^{\prime \prime} \, ) $, transformed respectively
by means of the matrix $T$ and its transpose $T^t$, satisfy the relations
 $ x^{\prime} y^{\prime} \,= \, q \, y^{\prime} x^{\prime}$ and
 $  x^{\prime \prime} y^{\prime \prime} \, = \, q \, y^{\prime \prime}
  x^{\prime \prime} \,  $ where
\begin{equation}
T  :  \left( \begin{array}{c} x \\ y \end{array} \right)\; \mapsto \;
  \left( \begin{array}{c} x^{\prime} \\ y^{\prime} \end{array} \right) \;
  = \left( \begin{array}{cc} A & B \\ C & D \end{array} \right)
  \left( \begin{array}{c} x \\ y \end{array} \right)   \label{eq:5}
\end{equation}
and
\begin{equation}
T^t :   \left( \begin{array}{c} x \\ y \end{array} \right)\; \mapsto \;
  \left( \begin{array}{c} x^{\prime \prime} \\ y^{\prime \prime} \end{array}
  \right) \;  = \left( \begin{array}{cc} A & C \\ B & D \end{array} \right)
  \left( \begin{array}{c} x \\ y \end{array} \right) \, .
	    \label{eq:6}
\end{equation}
What we emphasize here is that the relation in  Eq.(\ref{eq:2})  is
invariant not only under the transformation $T$  but also under its
transpose $T^t$ ( In this sense, a one-parameter quantum group
can be regarded as a
symmetry group of a quantum plane )
 and that {\em the generators of a quantum group are assumed
to be commutative with the coordinates of a quantum plane}.

In this work, we are naturally led to a two-parameter deformation of the 
group
$GL(2)$ and its corresponding quantum planes even though we do not
put any restriction
on the number of parameters at the outset.  Thus even though the
multi-parameter case has already been studied
\cite{ref14,ref17,ref18}, we shall concern ourselves with only the
two-parameter case in this work.
 Two-parameter quantum planes have still
attracted  attention  recently\cite{ref15,ref16}.

Now let us recall  two-parameter quantum groups.
 In fact, by solving the Yang-Baxter equation, one can get the
universal $R$-matrix
\begin{equation}
	      R_{p, \, q} \; = \;
		    \left( \begin{array}{cccc}
			   q & 0 & 0 & 0   \\
			   0 & 1 & 0 & 0   \\
      0 & q - \frac{1}{p} & \frac{q}{p} & 0 \\
			   0 & 0 & 0 & q
		    \end{array} \right)    \label{eq:7}
\end{equation}
where $p $ and $q$ are free  parameters\cite{ref9,ref10,ref11}.
From $RTT $ relations, one has the commutation relations
\begin{eqnarray}
  A \, B \; &=& \; p \, B \, A,   \hspace{2 cm}
			      C \, D \; = \; p \, D \, C,   \nonumber \\
  A \, C \; &=& \; q \, C \, A,   \hspace{2 cm}
			      B \, D \; = \; q \, D \, B,   \label{eq:8} \\
 p \, B \, C \; &=& q \, C \, B,  \hspace{2 cm}
      A \, D \, - \, D \, A \; = \; (p - \frac{1}{q} ) \, B \, C \, .
      \nonumber
\end{eqnarray}
We note that $R_{p, q } $ and   Eq. (\ref{eq:8})  become the
well-known $R_q $ solution and  Eq. (\ref{eq:1}),  respectively, in the
limit
 $p \, \rightarrow \, q$.
However,  Eq. (\ref{eq:2})  in the two-parameter case is not invariant
under
the transformation in Eq. (\ref{eq:6}). It is only invariant  under the
transformation in Eq. (\ref{eq:5}).  Whenever one requires that
 Eq. (\ref{eq:2})
be invariant under the two transformations  with the assumption that
the generators of a quantum group and the coordinates of a quantum plane
be commutative,
one is led to a one-parameter quantum group.

Our observation is that even though there are no restrictions
on the number of parameters at the outset,
one is led naturally to a two-parameter quantum group
    $GL_{p, q} (2)$
in such a manner that the commutation relations in Eq. (\ref{eq:8})
come directly from the condition that  $x \, y \; = \; q \, y \, x $
is preserved not only under the transformation in Eq. (\ref{eq:5})
but also under that in Eq.
(\ref{eq:6})  as in the one-parameter case, {\em if
one relaxes the commutation relations between
the generators of a quantum group  and the noncommuting
coordinates}.   Actually,
the remarkable fact is that
 even in the case of one-parameter quantum groups, the generators of a
quantum group  do not commute with the coordinates of the quantum plane
{\em generically}, as can be seen in the next section.

In Sec. II, we shall push
this observation further in a more general
fashion.   This formulation leads us to many probable quantum planes
associated with a quantum group.
  We shall discuss some special examples
in Sec. III.

\section{Two-parameter quantum group as a symmetry group }

Let
$ \left( \begin{array}{cl} A & B \\ C & D \end{array} \right) $
be an element of a  quantum group and
let us assume that for some numbers, $q_{ij}$'s,
\begin{eqnarray}
     x \, A \; &=& \; q_{11} \, A \, x,    \hspace{2 cm}
     y \, A \; = \; q_{21} \, A \, y,  \nonumber  \\
     x \, B \; &=& \; q_{12} \, B \, x,     \hspace{2 cm}
     y \, B \; = \; q_{22} \, B \, y, \nonumber \\
     x \, C \; &=& \; q_{13} \, C \, x,     \hspace{2 cm}
     y \, C \; = \; q_{23} \, C \, y,  \label{eq:12} \\
     x \, D \; &=& \; q_{14} \, D \, x,     \hspace{2 cm}
     y \, D \; = \; q_{24} \, D \, y \, .  \nonumber
\end{eqnarray}
Also let us assume that, under the transformations
in Eqs.(\ref{eq:5}) and (\ref{eq:6}), the
relation $x \, y \; = \; q \, y\, x$ is transformed, respectively, as
\begin{equation}
 x^{\prime} \, y^{\prime} \; = \; \bar{q} \, y^{\prime} \, x^{\prime}
      \label{eq:10}
\end{equation}
and
\begin{equation}
x^{\prime \prime} \, y^{\prime \prime} \;
      = \; \bar{\bar{q}} \, y^{\prime \prime} \, x^{\prime \prime} \, .
      \label{eq:11}
\end{equation}
Then, we have
\vspace{0.2cm}

 (1)
$ \left( \begin{array}{cl} A & B \\ C & D \end{array} \right)
      \in GL_{p,q^{\prime}}(2)$
 for some nonzero $p, q^{\prime} $ with $pq^{\prime} \ne -1 $,
\vspace{0.5cm}

(2) \hspace{1.5cm} $ \bar{q} = \bar{\bar{q}} $, \hspace{0.5cm} and
\vspace{0.2cm}

(3)
\vspace{-1.17cm}
\begin{eqnarray}
 q_{11} &=&   1,  \hspace{2.83cm}
 q_{21} = q\bar{q}^{-1} q_{14}  \; = \;q {q^{\prime}}^{-1} \, k, \nonumber \\
 q_{12} &=&  \bar{q} \, p^{-1}, \hspace{2.21cm}
	 q_{22} = q\bar{q} \, p^{-1} \,(\bar{q}
       - ( p - {q^{\prime}}^{-1} ) k ),       \label{eq:102}  \\
 q_{13} &=& \bar{q} \, {q^{\prime}}^{-1}, \hspace{2.21cm}
 q_{23} =  q\bar{q} \, {q^{\prime}}^{-1} \,
	       (\bar{q} - ( p - {q^{\prime}}^{-1} ) k ),
	  \nonumber \\
q_{14} &=& \bar{q} {q^{\prime}}^{-1} \, k, \hspace{2cm}
 q_{24} =  q{\bar{q}}^{2} \, {q^{\prime}}^{-1} p^{-1}
		(\bar{q} - ( p - {q^{\prime}}^{-1} ) k )
	 \nonumber
\end{eqnarray}
  where $k$ is a complex number. In this section, we shall prove the above
statement.
 The converse of the above statement is trivial.  Also we note that if one
requires that $q_{ij} = 1$, then $\bar{q} = p = q^{\prime} = q $
and
$ \left( \begin{array}{cl} A & B \\ C & D \end{array} \right)
      \in GL_{q}(2).$

The proof is as follows:
From the Eqs. (\ref{eq:10}) and (\ref{eq:11}), it follows that
\begin{eqnarray}
 & &  A \, C \; = \; q_1 \, C \, A,   \nonumber \\
 & &  B \, D \; = \; q_2 \, D \, B,       \label{eq:13} \\
     & &  q \, q_{14} \, A \, D \; - \; \bar{q} \, q_{21} \, D \, A \;
	= q\bar{q} \, q_{12} \, C \, B \;
	- \; q_{23} \, B \, C,           \nonumber
\end{eqnarray}
and
\begin{eqnarray}
 & &  A \, B \; = \; q_3 \, B \, A,   \nonumber \\
 & &  C \, D \; = \; q_4 \, D \, C,       \label{eq:14} \\
 & &  q \, q_{14} \, A \, D \; - \; \bar{\bar{q}} \, q_{21} \, D \, A \;
       = q\bar{\bar{q}} \, q_{13} \, B \, C \;
       - \; q_{22} \, C \, B \, , \nonumber
\end{eqnarray}
where $ q_1 \: = \: \bar{q} \, {q_{13}}^{-1} \, q_{11} , \;
 q_2 \: = \: \bar{q} \, {q_{24}}^{-1} \, q_{22} , \;
 q_3 \: = \: \bar{\bar{q}} \, {q_{12}}^{-1} \, q_{11},  $
and
     $ q_4 \: = \:  \bar{\bar{q}} \, {q_{24}}^{-1} \, q_{23} $.

We are now interested in those  $ q_{ij}$'s such that
     $ T \;=\; \left( \begin{array}{cc} A & B \\ C & D \end{array} \right) $
is an element of a quantum group. For the matrix $T$ to be such a matrix, the
entries  $ A, \, B, \,  C,$ and $D$ should be consistent with the
costructures of the Hopf algebra.
We note that the comultiplication $ \Delta $ and the antipode $ S $,
 among others,
satisfy the following relations:
\begin{eqnarray}
 & & \Delta \left( \begin{array}{cc}
	A & B \\
	C & D  \end{array} \right) \; = \;
	 \left( \begin{array}{cc}
	A & B \\
	C & D  \end{array} \right) \; \otimes \;
	  \left( \begin{array}{cc}
	A & B \\
	C & D  \end{array} \right)  \nonumber \\
    & &    = \;
      \left( \begin{array}{cc}
       A \otimes A \, + \, B \otimes C  & A \otimes B \, + \, B \otimes D \\
       C \otimes A \, + \, D \otimes C  & C \otimes B \, + \, D \otimes D
      \end{array} \right)  \label{eq:15}
\end{eqnarray}
and
\begin{equation}
     S \left( \begin{array}{cc}
	A & B \\
	C & D  \end{array} \right) \; = \;
	 \left( \begin{array}{cc}
	A & B \\
	C & D  \end{array} \right)^{-1}\, .  \label{eq:16}
\end{equation}
From the consistent conditions
   $\Delta \left( A  C \right) \; = \;
       q_1 \,\Delta \left( C  A \right)  $
and
   $\Delta \left( B  D \right) \; = \;
       q_2 \,\Delta \left( D  B \right) $,
we can have
    $q_1 \; = \; q_2 \; \equiv \;
    q^{\prime} $
and
\begin{equation}
    A \, D \; - \; D \, A \; = \; q^{\prime} \, C \, B \; - \;
     {q^{\prime} }^{-1} \, B \, C \, .    \label{eq:17}
\end{equation}
Also from the conditions
      $ \Delta \left( A  B \right) \; = \;
    q_3 \,\Delta \left( B A \right)  $
and
      $\Delta \left( C  D \right) \; = \;
      q_4 \,\Delta \left( D  C \right) $,
it follows that
    $q_3 \; = \; q_4 \; \equiv \; p $
and
\begin{equation}
     A \, D \; - \; D \, A \; = \; p \, B \, C \; - \;
     p^{-1} \, C \, B \, .    \label{eq:18}
\end{equation}
From  Eqs. (\ref{eq:17}) and (\ref{eq:18}), it follows that
\begin{equation}
 p \, B \, C \; = \; q^{\prime} \, C \, B \, ,  \label{eq:20}
\end{equation}
unless $ pq^{\prime} = -1 $.
Thus, we construct a two-parameter deformation of $GL(2)$:
\begin{eqnarray}
  A \, B \; &=& \; p \, B \, A,   \hspace{2 cm}
			      C \, D \; = \; p \, D \, C,   \nonumber \\
  A \, C \; &=& \; q^{\prime} \, C \, A,   \hspace{2 cm}
       B \, D \; = \; q^{\prime} \, D \, B,   \label{eq:100} \\
 p \, B \, C \; &=& q^{\prime} \, C \, B,  \hspace{2 cm}
       A \, D \, - \, D \, A \; = \; (p - \frac{1}{q^{\prime}} ) \,
       B \, C \, .
       \nonumber
\end{eqnarray}
Hence 
$
\left( \begin{array}{cc}
	A & B \\
	C & D  \end{array} \right) \in GL_{p,q^{\prime}} $.

Next, Eq. (\ref{eq:16}) implies the existence of the inverse matrix
$ T^{-1} $.
From the ansatz
\begin{equation}
    \left( \begin{array}{cc}
    A & B \\
    C & D \end{array} \right) \,
    \left( \begin{array}{cc}
    D & \beta \, B \\
    \gamma \, C & \alpha  \, A \end{array} \right) \, {\cal D}^{-1} \; = \;
   \left( \begin{array}{cc}
   1 & 0 \\
    0 & 1   \end{array} \right)\, ,   \nonumber
\end{equation}
we can find
    $ \alpha \, = 1 , \: \: \beta \, =  - p^{-1}  ,
    \:\: \gamma \, =  - \, p, $ and
    $ {\cal D} \, = \, A D \, - \, p  B  C \, = \, D A \, - \, p^{-1}
    C B \, , $ which is consistent with Eq. (\ref{eq:18}).

The quantum determinant ${\cal D} $ satisfies
\begin{eqnarray}
A{\cal D}& =& {\cal D}A,  \hspace{3cm}  B{\cal D}  
		    = p^{-1}q^{\prime}{\cal D}B,
	       \nonumber  \\
C{\cal D}& =& p{q^{\prime}}^{-1}{\cal D}C, 
	   \hspace{2cm}  D{\cal D}  = {\cal D}D. \label{eq:1001} 
\end{eqnarray}
This gives us 
\begin{eqnarray}
   {\cal D}^{-1} A \;& = & \; A {\cal D}^{-1},   \nonumber \\
   {\cal D}^{-1} B \;& = & \; q^{\prime} p^{-1} B  {\cal D}^{-1}, \nonumber \\
   {\cal D}^{-1} C \;& = & \; p {q^{\prime}}^{-1}  C  {\cal D}^{-1},  
		\label{eq:19} \\
   {\cal D}^{-1} D \;& = & \; D  {\cal D}^{-1}, \nonumber
\end{eqnarray}
which is consistent with
the requirement
\begin{equation} \left( \begin{array}{cc}
	D & - \frac{1}{p} B \\
       -p \, C & A \end{array} \right) \; {\cal D}^{-1} \;
	      \left( \begin{array}{cc}
	A & B \\
	C & D \end{array} \right) \; = \;
       \left( \begin{array}{cc}
	1 & 0 \\
	0 & 1 \end{array} \right) .
\end{equation}
The result in Eq. (\ref{eq:19}) is the same as the one in Ref 12.

Furthermore, the third
equations in Eqs. (\ref{eq:13}) and (\ref{eq:14})  also should be
identical to  Eq. (\ref{eq:18}).
If $qq_{14} \ne \bar{q}q_{21} $, the third eq. in Eq. (\ref{eq:13})
is $qq_{14}AD - \bar{q}q_{21}DA = \xi BC $ 
where $\xi = q\bar{q}q_{12}{q^{\prime}}^{-1}p - q_{23} $. In the case when
$\xi = 0 $, $qq_{14}AD = \bar{q}q_{21}DA $ which is of the form
$AD = \epsilon DA $ with $\epsilon \ne 1 $. 
However, the relation $\Delta (AD) = \epsilon \Delta (DA)$ 
leads us to $\epsilon = 1 $, which is a contradiction.
When $\xi \ne 0 $, we have two cases: $p \ne {q^{\prime}}^{-1} $  
and $ p = {q^{\prime}}^{-1} $. The first case together with Eq. (\ref{eq:18})
gives $(qq_{14}(p -{q^{\prime}}^{-1}) - \xi)AD = 
(\bar{q}q_{21}(p-{q^{\prime}}^{-1})-\xi)DA $.
In every possible case, this equation
  contradicts  either the fact that ${\cal D} = AD - pBC$ is
 invertible or that
 $\epsilon = 1 $ from
 $\Delta (AD) = \epsilon \Delta (DA)$ as in the above.
 In the second case when $ p = {q^{\prime}}^{-1}$, 
 $AD = DA =  \delta BC $
for some number $\delta $.  However,
from the relation $\Delta(AD) = \delta \Delta(BC)$, $\delta = p$, which is 
a contradiction to the existence of ${\cal D}$.  Thus, we conclude that
$qq_{14} = \bar{q}q_{21}$.   The equation $qq_{14} = \bar{\bar{q}}q_{21} $
follows from the third equation in Eq. (\ref{eq:14}) by a completely
analogous method.   Hence, we have
\begin{equation}
\bar{q} = \bar{\bar{q}}. \label{eq:1002} 
\end{equation}

Now let us summarize the 
 relations between the $q_{ij}$'s:
\begin{eqnarray}
   & q^{\prime} \;& =  \; \bar{q} \, {q_{13}}^{-1} \, q_{11}
	\; = \; \bar{q} \, {q_{24}}^{-1} \, q_{22},  \nonumber \\
     &  p \;& = \;\bar{q} \, {q_{12}}^{-1} \, q_{11}
	\; = \; \bar{q} \, {q_{24}}^{-1} \, q_{23},  \nonumber \\
   & q&q_{14}  \; = \; \bar{q}q_{21},     \label{eq:101} \\
 &   p &\; - \; {q^{\prime}}^{-1} \; =\;
 \bar{q} \, {q_{14}}^{-1} \, q_{13} \; - \; q^{-1} \, {q^{\prime}}^{-1} \,
p \; {q_{14}}^{-1} \, q_{22} . \nonumber
\end{eqnarray}
The relation between $q, \bar{q}$, and $q^{\prime}$ depends
 on the choice of $q_{ij}$'s.
There may be (infinitely) many choices for $q_{ij}$'s consistent with the
theory of  quantum groups.
In effect, there are two unknowns  since there are
six independent relationships between them, as can be seen in Eq.
(\ref{eq:101}).
Without loss of generality, we may assume that $q_{14} \; = \; k q_{13}$ for
some number $k$. Then we can express all of the $q_{ij}$'s in one unknown
$q_{11}$, which may be regarded as a proportional constant.  Hence,
if we put $q_{11} \; = \; 1 $ for simplicity, we have

 \begin{eqnarray}
 q_{11} &=& 1,
  \hspace{2.83cm}
 q_{21} = q{\bar{q}}^{-1} q_{14}  \; = 
	     \;q {q^{\prime}}^{-1} \, k, \nonumber \\
q_{12} & =& \bar{q} \, p^{-1},
 \hspace{2.2 cm}
      q_{22} = q\bar{q} \, p^{-1} \,(\bar{q} - ( p - {q^{\prime}}^{-1} ) k ),
      \label{eq:1003}  \\
 q_{13} & = &\bar{q} \, {q^{\prime}}^{-1},
  \hspace{2.2 cm}
   q_{23} =  q\bar{q} \, {q^{\prime}}^{-1} \,
	       (\bar{q} - ( p - {q^{\prime}}^{-1} ) k ),
	  \nonumber \\
 q_{14} &=& \bar{q} {q^{\prime}}^{-1} \, k,
\hspace{2 cm} 
q_{24} =  q{\bar{q}}^{2} \, {q^{\prime}}^{-1} p^{-1} 
		(\bar{q} - ( p - {q^{\prime}}^{-1} ) k ) 
	 \nonumber
\end{eqnarray}
where $k$ is the only parameter to be determined.
Thus, we prove the statement.  As seen in the above, the choice $q_{11} = 1$
is arbitrary.  In other words, the assumption that the generators 
of a one-parameter quantum group commute with the coordinates of the quantum
plane is very special.   They do not commute generically.  

From  Eq. (\ref{eq:100}), it is obvious that
  $\left( \begin{array}{cc}  A & B  \\ C & D \end{array} \right)
  \in GL_{p,q^{\prime}} (2) $
if and only if
  $\left( \begin{array}{cc}  A & C  \\ B & D \end{array} \right)
  \in GL_{q^{\prime},p} (2) $.
On the other hand, 
$GL_{p,q^{\prime}} (2) \; = \; GL_{q^{\prime},p} (2)$ in the sense that
$GL_{p,q^{\prime}} (2) $ 
and $ GL_{q^{\prime},p} (2)$ are the algebras freely generated by
$A, \; B, \; C, \; D$, and ${\cal D}^{-1} $ modulo the relations given by
 Eqs. (\ref{eq:100}) and (\ref{eq:19}) and by the equations
  $(AB \; - \; pBC) {\cal D}^{-1} \; - \; 1  $ and
  $      {\cal D}^{-1} (AB \; - \; p BC) \; - \; 1 $.
Thus,  Manin's viewpoint that quantum groups are symmetry groups of
quantum planes is recovered as in the one-parameter case under the
commutation relation in Eq. (\ref{eq:12}) with $q_{ij}$'s  given by 
Eq. (\ref{eq:1003}) between quantum group generators and noncommuting
coordinates.

\section{Quantum Planes Associated With A Quantum Group }

In this section, we shall  
discuss several interesting choices of $q_{ij}$'s.
The diversity of the choices of  $q_{ij}$'s means the diversity of 
quantum planes for a given quantum group.
\vspace{1cm}
 
{\bf Case I: } $\bar{q} = q $
 
This case corresponds to the standard way of dealing with 
quantum planes.  
Then,
\begin{eqnarray}
 q_{11}& =& 1,
 \hspace{2.83 cm}
 q_{21} =  q_{14}  \; = \;q {q^{\prime}}^{-1} \, k, \nonumber \\
q_{12} &=& q \, p^{-1},
  \hspace{2.2 cm}
      q_{22} = q^2 \, p^{-1} \,(q - ( p - {q^{\prime}}^{-1} ) k ),
      \label{eq:1004}  \\
 q_{13} &=& q \, {q^{\prime}}^{-1},
 \hspace{2.2 cm}
   q_{23} =  q^2 \, {q^{\prime}}^{-1} \,
	       (q - ( p - {q^{\prime}}^{-1} ) k ),
	  \nonumber \\
 q_{14}& =& q {q^{\prime}}^{-1} \, k,
\hspace{2 cm} 
q_{24} =  q^3 \, {q^{\prime}}^{-1} p^{-1} 
		(q - ( p - {q^{\prime}}^{-1} ) k ) \, .
	 \nonumber
\end{eqnarray}
 Now we introduce the  exterior differential $d$ 
as in Ref 15 and 16
  except for the following:
\begin{eqnarray}
       (dx) \, A \; = \; q_{11} \, A \, dx,
    & \hspace{2 cm} &
       (dy) \, A \; = \; q_{21} \, A \, dy, \nonumber \\
       (dx) \, B \; = \; q_{12} \, B \, dx,
    & \hspace{2 cm} &
       (dy) \, B \; = \; q_{22} \, B \, dy, \nonumber \\
   (dx) \, C \; = \; q_{13} \, C \, dx,
    & \hspace{2 cm} &
       (dy) \, C \; = \; q_{23} \, C \, dy, \label{eq:26}  \\
       (dx) \, D \; = \; q_{14} \, D \, dx,
    & \hspace{2 cm} &
       (dy) \, D \; = \; q_{24} \, D \, dy  \nonumber
\end{eqnarray}
where
    $ q_{ij}$'s
 satisfy Eq. (\ref{eq:1004}).

Now if we require that $dx \, dy \; = \; - \frac{1}{p} \, dy \, dx$ is
preserved under the transformation $T$, it is easy to see that
$k \; = \; \frac{ q^{\prime} ( q p - 1)}{p(q^{\prime} p -1)} $.
Thus, we have, with $q_{11} = 1$,
\begin{eqnarray}
q_{12}& = q \, p^{-1},
    \hspace{2 cm}&
   q_{21} \; = \; q_{14}  \; = 
	  \;\frac{ q ( q p - 1)}{p(q^{\prime} p -1)}, \nonumber \\
   q_{13} &= q \, {q^{\prime}}^{-1},  \hspace{2 cm}&
   q_{22} = q^2 \, p^{-2},   \label{eq:1005}  \\
   q_{14} & = \frac{ q ( q p - 1)}{p(q^{\prime} p -1)},
    \hspace{2 cm}&
   q_{23} =  q^2 \, {q^{\prime}}^{-1} \,p^{-1},    \nonumber \\
   &\hspace{2 cm} &
q_{24} =  q^3 \, {q^{\prime}}^{-1} p^{-2}  \,.
    \nonumber
\end{eqnarray}
Now we may go further. In fact, it is natural to require that the
two-parameter case become the one-parameter case in some limit.
Therefore, if $q_{ij} \longrightarrow 1 $ as $p  \longrightarrow q^{\prime}$,
then we must set $q^{\prime} = q$.
Hence, Eq. (\ref{eq:100}) is the same as Eq. (\ref{eq:8}), and Eq.
(\ref{eq:1005}) becomes
\begin{eqnarray}
& & q_{11} \; = \; q_{13} \; = \; 1, \nonumber \\
& & q_{12} \; = \; q_{14} \; = \; q_{21} \; = \; q_{23} \; = \;
   q \, p^{-1},    \label{eq:25}  \\
& & q_{22} \; = \; q_{24} \; = \; q^2 \, p^{-2} \, .  \nonumber
\end{eqnarray}

 The virtue of this formulation is that the relations for the differentials
on a quantum plane are preserved not only under $T$ but also under $T^t$.
According to Ref 16,
one can define the differential calculus on a quantum plane in the
one-parameter case:
For an exterior differential $d$ which is linear and satisfies $d^2 = 0$
and the
Leibnitz rule, one can choose
\begin{eqnarray}
 dx \, dy \; &=& \; - \frac{1}{q} \, dy \, dx,   \nonumber \\
	 x \, dx \; &=& \; q^2 \, (dx) \, x,         \nonumber \\
	 x \, dy \; &=& \; q \,
	  (dy) \, x \; + \;  ( q^2 -1 ) \, (dx) \, y,   \label{eq:4}      \\
	y \, dx \; &=& \; q \, (dx) \, y,          \nonumber \\
 y \, dy \; &=& \; q^2 \, (dy) \, y \, \, .   \nonumber
\end{eqnarray}
Also, by the same method as in the one-parameter case, we obtain the
following relations for the differentials in the two-parameter case:

\begin{eqnarray}
 dx \, dy \; &=& \; - \frac{1}{p} \, dy \, dx,   \nonumber \\
	 x \, dx \; &=& \; p \, q \, (dx) \, x,         \nonumber \\
       x \, dy \; &=& 
	      \; q \, (dy) \, x \; + \;  ( p \, q -1 ) \, (dx) \, y,
					\label{eq:9}      \\
	y \, dx \; &=& \; p \, (dx) \, y,          \nonumber \\
 y \, dy \; &=& \; p \, q \, (dy) \, y \, \, .   \nonumber
\end{eqnarray}
We note that Eq. (\ref{eq:4}) is invariant under the transformations
 $ T  $ and  $T^t $.
Eq. (\ref{eq:9}) is also invariant under the transformation
  $T $,  but it is easy to see that it is not invariant under the
transformation
  $T^t$
if the quantum group generators $ A\, , \; B \, , \; C $, and $D$ commute
with the noncommuting coordinates $ x\, , \: y \, .  $
However, if we choose the $q_{ij}$'s as in Eq. (\ref{eq:25}),
then a lengthy but straightforward calculation shows the nice property that
Eq. (\ref{eq:9}) is invariant not only under the transformation
 $T$
but  also under the transformation
  $T^t$.
Moreover, we have $dx^{\prime} dx^{\prime} \; = \; dy^{\prime} dy^{\prime}
\; = \; 0$ and $dx^{\prime \prime} dx^{\prime \prime} \; = \;
dy^{\prime \prime} dy^{\prime \prime} \; = \; 0 $.
\vspace{1cm}
 
{\bf Case II: }  $ p = q^{\prime} $ 
\vspace{0.5cm}

Let   $\left( \begin{array}{cc}  A & B  \\ C & D \end{array} \right)
  \in GL_{q^{\prime}} (2) $. 
If we put $ k = q_{12}$ ( In effect, this choice of $k$ gives the same
equation, Eq. (\ref{eq:25}), as case I  above  ), then
  \begin{eqnarray}
 q_{11} &= & 1,
  \hspace{2.83 cm}
 q_{21} = q\bar{q}{q^{\prime}}^{-1}p^{-1}, \nonumber \\
q_{12} & = & \bar{q} \, p^{-1},
  \hspace{2.2cm}
      q_{22} = q\bar{q}^{2}{q^{\prime}}^{-1}p^{-2}, 
      \label{eq:1010}  \\
 q_{13} & = & \bar{q} {q^{\prime}}^{-1},
  \hspace{2.2 cm}
   q_{23} =  q\bar{q}^{2} \, {q^{\prime}}^{-2}p^{-1},
	  \nonumber \\
 q_{14} &= &\bar{q}^{2} {q^{\prime}}^{-1}p^{-1},
 \hspace{1.5 cm} 
q_{24} =  q{\bar{q}}^{3}  {q^{\prime}}^{-2}p^{-2} 
	 \nonumber
\end{eqnarray}
In order to see interesting aspects of quantum planes, it is enough only to
 consider the one-parameter case.
Thus,
 if put $p = q^{\prime}$,
 \begin{eqnarray}
 q_{11}& = & 1,
  \hspace{2.83cm}
 q_{21} = q\bar{q}{q^{\prime}}^{-2}, \nonumber \\
q_{12} & = & \bar{q} \, {q^{\prime}}^{-1},
  \hspace{2.2 cm}
      q_{22} = q\bar{q}^{2}{q^{\prime}}^{-3}, 
      \label{eq:1011}  \\
 q_{13} &=& \bar{q} {q^{\prime}}^{-1},
  \hspace{2.2 cm}
   q_{23} =  q\bar{q}^{2} {q^{\prime}}^{-3},
	  \nonumber \\
 q_{14} &=& \bar{q}^{2} {q^{\prime}}^{-2},
 \hspace{2 cm} 
q_{24} =  q{\bar{q}}^{3}  {q^{\prime}}^{-4}. 
	 \nonumber
\end{eqnarray}
The quantum plane such that $xy = qyx $
corresponding to these values of the $q_{ij}$'s
is transformed
 into $x^{\prime}y^{\prime} = \bar{q}y^{\prime }x^{\prime} $ and
 $x^{\prime \prime}y^{\prime \prime}
   = \bar{q}y^{\prime \prime}x^{\prime \prime} $, 
respectively, under the action  
 $\left( \begin{array}{cc}  A & B  \\ C & D \end{array} \right) $ 
and its tranpose.

 Now if we take  a quantum plane for $GL_{q^{\prime}}$ such that
$ q = 1 $ and $\bar{q} = q^{\prime}$,
then
\begin{equation}
 q_{1i} = 1,  \hspace{2cm} q_{2i} = {q^{\prime}}^{-1}, 
\end{equation}
for $i = 1, \cdots, 4 $.
This quantum plane is generated by $ x, y $ such that $xy = yx $ and is 
transformed as $ x^{\prime}y^{\prime} = q^{\prime}y^{\prime}x^{\prime}$.
However, $x^{\prime}, y^{\prime} $ do not obey Eq. (\ref{eq:12}).
The case when $q = \bar{q} = 1$ is also interesting since 
the quantum plane looks like an ordinary plane in the sense that
it is  generated by commuting coordinates. 
If we take  a  quantum plane for $GL_{q^{\prime}}$ such that
$ q = q^{\prime} $ and $\bar{q} = q^{\prime}$,
then $q_{ij} = 1 $.  This quantum plane is the original one
 \cite{ref1}.  

\section{Conclusions}

In the one-parameter case, the condition that
  $x \, y \; = \; q \, y \, x $
is preserved under the transformation $T$
and its transpose $T^t$ gives 
the commutation relation between the generators
of a quantum group $GL_{q} (2) $.
Here, one assumes that
 the generators of a quantum group  commute with the noncommuting
coordinates of a quantum plane.

In this work, we have relaxed the assumption and investigated its
consequences.
We are naturally led to a two-parameter deformation of the group
$GL(2)$ and its corresponding quantum planes even though we do not
put any restrictions at the outset 
on the number of parameters.
As a by-product, this formulation supports
 Manin's viewpoint that quantum groups are symmetry groups of quantum
planes,
and the diversity of the choices of $q_{ij}$'s 
shows that there can be many quantum planes for a given
quantum group $GL_{p,q}$.
Associated with a given quantum group,
there are some special quantum planes such as the original one
 in the literature.
Especially, a quantum plane which looks like an ordinary plane
attracts much attention and seems to be worthy of further research.

{\bf ACKNOWLEDGMENTS}

This work was supported by Ministry of Education, Project
No. BSRI-95-2414, a grant from TGRC-KOSEF 1994 and Korea Research
Center for Theoretical Physics and Chemistry. 

\end{document}